\documentclass{llncs}
\usepackage{makeidx}  
\usepackage{amsmath}
\usepackage{amssymb}
\usepackage{amsfonts}
\usepackage[english]{babel}
\usepackage{fancyhdr}

\def \egal {\stackrel{{\rm def}}{=}}
\def\tf{\tilde{f}_n}
\def\hf{\hat{f}_n}

\newcommand \cA{{\cal A}}

\newcommand \cC{{\cal C}}

\newcommand \cF{{\cal F}}
\newcommand \cG{{\cal G}}

\newcommand \cN{{\cal N}}

\newcommand \cP{{\cal P}}

\newcommand \cX{{\cal X}}

\newcommand{\1}{{\rm 1}\kern-0.24em{\rm I}}

\begin{document}

\title{Optimal oracle inequality for aggregation of classifiers under low noise condition}
\titlerunning{Aggregation Method in Classification}  
%
\author{Guillaume Lecu\'e \footnote{Paper to be considered for the Mark Fulk Award for the "best student paper".}}
\authorrunning{Guillaume Lecu\'e}   
%
\tocauthor{Guillaume Lecu\'e (Universit\'e Paris VI)}
\institute{Laboratoire de Probabilit\'es et Mod\`eles Al\'eatoires
(UMR CNRS 7599)\\Universit\'e Paris VI\\ 4 pl.Jussieu, BP 188,
75252 Paris, France ,\\
\email{lecue@ccr.jussieu.fr}}

\maketitle              

\begin{abstract} \noindent
We consider the problem of optimality, in a minimax sense, and
adaptivity to the margin and to regularity in binary
classification. We prove an oracle inequality, under the margin
 assumption (low noise condition), satisfied by an aggregation
procedure which uses exponential weights. This oracle inequality
has an optimal residual: $(\log M/n)^{\kappa/(2\kappa-1)}$ where
$\kappa$ is the margin parameter, $M$ the number of classifiers to
aggregate and $n$ the number of observations. We use this
inequality first to construct minimax classifiers under margin and
regularity assumptions and second to aggregate them to obtain a
classifier which is adaptive both to the margin and regularity.
Moreover, by aggregating plug-in classifiers (only $\log n$), we
provide an easily implementable classifier adaptive both to the
margin and to regularity.
\end{abstract}

\section{Introduction}

Let $(\cX,\cA)$ be a measurable space. We consider a random
variable $(X,Y)$ with values in $\cX\times\{-1,1\}$ and denote by
$\pi$ the distribution of $(X,Y)$. We denote by $P^X$ the marginal
of $\pi$ on $\cX$ and $\eta(x)=\mathbb{P}(Y=1|X=x)$ the
conditional probability function of $Y=1$ given that $X=x$. We
denote by $D_n=(X_i,Y_i)_{i=1,\ldots,n}$, $n$ i.i.d. observations
of the couple $(X,Y)$.

We recall some usual notions introduced for the classification
framework. A {\it{prediction rule}} is a measurable function
$f:\cX\longmapsto\{-1,1\}$. The {\it{misclassification error}}
associated to $f$ is $$R(f) = \mathbb{P}(Y\neq f(X)).$$ It is well
known (see, e.g., \cite{dgl:96}) that $\min_f R(f) = R(f^*)\egal
R^{*},$ where the prediction rule $f^*$ is called {\it Bayes rule}
and is defined by $$ f^*(x)= {\rm{sign}}(2\eta(x)-1).$$The minimal
risk $R^*$ is called the {\it{Bayes risk}}. A {\it{classifier}} is
a  function, $\hat{f}_n=\hat{f}_n(X,D_n)$, measurable with respect
to $D_n$ and $X$ with values in $\{-1,1\}$, that assigns to the
sample $D_n$ a prediction rule
$\hat{f}_n(.,D_n):\cX\longmapsto\{-1,1\}$. A key characteristic of
$\hf$ is the value of {\it generalization error}
$\mathbb{E}[R(\hf)]$. Here$$R(\hf) = \mathbb{P}(Y\neq \hf(X)
|D_n).$$ The performance of a classifier $\hf$ is measured by the
value $\mathbb{E}[R(\hf) - R^*]$ called the {\it{excess risk}} of
$\hf$. We say that the classifier $\hf$ learns with the
convergence rate $\phi(n)$, where $(\phi(n))_{n\in\mathbb{N}}$ is
a decreasing sequence, if there exists an absolute constant $C>0$
such that for any integer $n$, $\mathbb{E}[R(\hf) - R^*]\leq
C\phi(n)$. Theorem 7.2 of \cite{dgl:96} shows that no classifier
can learn with a given convergence rate for arbitrary underlying
probability distribution $\pi$.

In this paper we focus on entropy assumptions which allow us to
work with finite sieves. Hence, we first work with a finite model
for $f^*$: it means that we take a finite class of prediction
rules $\cF=\{f_1,\ldots,f_M\}$. Our aim is to construct a
classifier $\hf$ which mimics the best one of them w.r.t. to the
excess risk and with an optimal residual. Namely, we want to state
an oracle inequality
\begin{equation}\label{exampleoracleinequality}\mathbb{E}\left[R(\hf)-R^* \right]\leq
a_0\min_{f\in\cF}(R(f)-R^*)+C\gamma(M,n),\end{equation} where
$a_0\geq1$ and $C>0$ are some absolute constants and $\gamma(M,n)$
is the residual. The classical procedure, due to Vapnik and
Chervonenkis (see, e.g. \cite{dgl:96}), is to look for an ERM
classifier,i.e., the one which minimizes the {\it{empirical risk}}
\begin{equation}\label{ER}
R_n(f)=\frac{1}{n}\sum_{i=1}^n\1_{\left\{Y_if(X_i)\leq0\right\}},
\end{equation}
over all prediction rules $f$ in $\cF$, where $\1_E$ denotes the
indicator of the set $E$. This procedure leads to optimal
theoretical results (see, e.g. Chapter 12 of \cite{dgl:96}), but
minimizing the empirical risk (\ref{ER}) is computationally
intractable for sets $\cF$ of classifiers with large cardinality
(often depending on the sample size $n$), because this risk is
neither convex nor continuous. Nevertheless, we might base a
tractable estimation procedure on minimization of a convex
surrogate $\phi$ for the loss ( \cite{lv:04}, \cite{by:02},
\cite{bbl:05}, \cite{blv:03}, \cite{ss:04} and \cite{ss:05}). A
wide variety of classification methods in machine learning are
based on this idea, in particular, on using the convex loss
associated to support vector machines (\cite{cv:95},
\cite{ssbook:02}),$$\phi(x)=\max(0,1-x),$$ called the {\it
hinge-loss}. The risk associated to this loss is called the {\it
hinge risk} and is defined by
$$A(f)=\mathbb{E}[\max(0,1-Yf(X))],$$ for all
$f:\cX\longmapsto\mathbb{R}$. The {\it{optimal hinge risk}} is
defined by
\begin{equation}\label{OHR}
A^*=\inf_f A(f),
\end{equation}
where the infimum is taken over all measurable functions $f$. The
Bayes rule $f^*$ attains the infimum in (\ref{OHR}) and, moreover,
denoting by $R(f)$ the misclassification error of ${\rm sign}(f)$
for all measurable functions $f$ with values in $\mathbb{R}$,
Zhang, cf. \cite{z:04}, has shown that,\begin{equation}\label{e1}
 R(f)-R^*\leq A(f)-A^*,\end{equation}for any real valued measurable function $f$. Thus, minimization of
the {\it{excess hinge risk}} $A(f)-A^*$ provides a reasonable
alternative for minimization of the excess risk. In this paper we
provide a procedure which does not need any minimization step. We
use a convex combination of the given prediction rules, as
explained in section \ref{oraclesection}.

The difficulty of classification is closely related to the
behavior of the conditional probability function $\eta$ near $1/2$
(the random variable $|\eta(X)-1/2|$ is sometimes called the
theoretical margin). Tsybakov  has introduced, in \cite{tsy:04},
an assumption on the the margin, called {\it{margin (or low noise)
assumption}},

 \noindent{\it{{\bf{(MA) Margin (or low noise) assumption.}}
 The probability distribution $\pi$ on
the space $\cX \times \{-1,1 \}$ satisfies the margin assumption
MA($\kappa$) with margin parameter $1\leq\kappa<+\infty$ if there
exists $c_0>0$ such that,
 \begin{equation}\label{e4}
   \mathbb{E}\left\{ |f(X)-f^{*}(X)|
\right\}\leq c_0 \left(R(f)-R^{*} \right)^{1/\kappa},
\end{equation}
for all measurable functions $f$ with values in $\{-1,1\}.$ }}

\noindent Under this assumption, the risk of an ERM classifier
over some fixed class $\cF$ can converge to the minimum risk over
the class with {\it{fast rates}}, namely faster than $n^{-1/2}$
(cf. \cite{tsy:04}). On the other hand, with no margin assumption
on the joint distribution $\pi$ (but combinatorial or complexity
assumption on the class $\cF$), the convergence rate of the excess
risk is not faster than $n^{-1/2}$ (cf. \cite{dgl:96}).

\noindent In this paper we suggest an easily implementable
procedure of aggregation of classifiers and prove the following
results:
 \begin{enumerate}
    \item We obtain an oracle inequality for our procedure and we use
it to show that our classifiers are adaptive both to the margin
parameter (low noise exponent) and to a complexity parameter.
    \item We generalize the lower bound inequality stated in Chapter 14 of
\cite{dgl:96}, by introducing the margin assumption and deduce
optimal rates of aggregation under low noise assumption in the
spirit of Tsybakov \cite{tsy:03}.
    \item We obtain classifiers with minimax
fast rates of convergence on a H\"older class of conditional
probability functions $\eta$ and under the margin assumption.
 \end{enumerate}

The paper is organized as follows. In Section $2$ we prove an
oracle inequality for our convex aggregate, with an optimal
residual, which will be used in Section \ref{example} to construct
minimax classifiers and to obtain adaptive classifiers by
aggregation of them. Proofs are given in Section \ref{proofs}.

\par

\section{Oracle Inequality}\label{oraclesection}

We have $M$ prediction rules $f_1,\ldots,f_M$. We want to mimic
the best of them according to the excess risk under the margin
assumption. Our procedure is using exponential weights. Similar
constructions in other context can be found, e.g., in
\cite{abb:97}, \cite{y:00}, \cite{catbook:01}, \cite{bl:04},
\cite{lec:05}, \cite{lec:06}, \cite{v:90}. Consider the following
aggregate which is a convex combination with exponential weights
of $M$ classifiers,
\begin{equation}\label{aggregate}
 \tilde{f_n}=\sum_{j=1}^{M}w_j^{(n)} f_j,
\end{equation}
where
 \begin{equation}\label{coef1}
w_j^{(n)}=\frac{\exp\left(
\sum_{i=1}^nY_if_j(X_i)\right)}{\sum_{k=1}^M\exp\left(
\sum_{i=1}^nY_if_k(X_i)\right)},\quad \forall j=1,\ldots,M.
 \end{equation}
Since $f_1,\ldots,f_M$ take their values in $\{-1,1\}$, we have,
 \begin{equation}\label{coefrisk}
 w_j^{(n)}=\frac{\exp\left( -nA_n(f_j)\right)}{\sum_{k=1}^{M}\exp\left(
 -nA_n(f_k)\right)},
 \end{equation}
for all $j\in\{1,\ldots,M\}$, where
 \begin{equation}\label{EHR}
A_n(f)=\frac{1}{n}\sum_{i=1}^n\max(0,1-Y_if(X_i))
\end{equation} is the
empirical analog of the hinge risk. Since $A_n(f_j)=2R_n(f_j)$ for
all $j=1,\ldots,M$, these weights can be written in terms of the
empirical risks of $f_j$'s,
$$w_j^{(n)}=\frac{\exp\left( -2nR_n(f_j)\right)}{\sum_{k=1}^{M}\exp\left(
-2nR_n(f_k)\right)}, \ \forall j=1,\ldots,M.$$ Remark that, using
the definition (\ref{coefrisk}) for the weights, we can aggregate
functions with values in $\mathbb{R}$ (like in theorem
\ref{theooracleA}) and not only functions with values in
$\{-1,1\}$.

The aggregation procedure defined by (\ref{aggregate}) with
weights (\ref{coefrisk}), that we can called aggregation with
exponential weights (AEW), can be compared to the ERM one. First,
our AEW method does not need any minimization algorithm contrarily
to the ERM procedure. Second, the AEW is less sensitive to the
over fitting problem. Intuitively, if the classifier with smallest
empirical risk is over fitted (it means that the classifier fits
too much to the observations) then the ERM procedure will be over
fitted. But, if other classifiers in $\cF$ are good classifiers,
our procedure will consider their "opinions" in the final decision
procedure and these opinions can balance with the opinion of the
over fitted classifier in $\cF$ which can be false because of its
over fitting property. The ERM only considers the "opinion" of the
classifier with the smallest risk, whereas the AEW takes into
account all the opinions of the classifiers in the set $\cF$. The
AEW is more temperate contrarily to the ERM. Understanding why
aggregation procedure are often more efficient than the ERM
procedure from a theoretical point of view is a deep question, on
which we are still working at this time this paper is written.
Finally, the following proposition shows that the AEW has similar
theoretical property as the ERM procedure up to the residual
$(\log M)/n$.
\begin{proposition}\label{propAEWERM}
Let $M\geq2$ be an integer, $f_1,\ldots,f_M$ be $M$ real valued
functions on $\cX$. For any integers $n$, the aggregate defined in
(\ref{aggregate}) with weights (\ref{coefrisk}) $\tilde{f_n}$
satisfies
  $$A_n(\tilde{f_n})\leq
  \min_{i=1,\ldots,M}A_n(f_i)+\frac{\log(M)}{n}.$$
\end{proposition}
The following theorem provides first an exact oracle inequality
w.r.t. the hinge risk satisfied by the AEW procedure and second
shows its optimality among all aggregation procedures. We deduce
from it that, for a margin parameter $\kappa\geq1$ and a set of
$M$ functions with values in $[-1,1]$, $\cF=\{f_1,\ldots,f_M\}$,
$$\gamma(\cF,\pi,n,\kappa)=\sqrt{\frac{\min_{f\in\cF}(A(f)-A^*)^{\frac{1}{\kappa}}\log
M}{n}}+\left(\frac{\log M}{n}\right)^{\frac{\kappa}{2\kappa-1}}$$
is an optimal rate of convex aggregation of $M$ functions with
values in $[-1,1]$ w.r.t. the hinge risk, in the sense of
\cite{lec:06}.
\begin{theorem}[{\bf Oracle inequality and Lower bound}]\label{theooracleA}
Let $\kappa\geq1$. We assume that $\pi$ satisfies MA($\kappa$). We
denote by $\cC$ the convex hull of a finite set of functions with
values in $[-1,1]$, $\cF=\{f_1,\ldots,f_M\}$. The AEW procedure,
introduced in (\ref{aggregate}) with weights (\ref{coefrisk})
(remark that the form of the weights in (\ref{coefrisk}) allows to
take real valued functions for the $f_j$'s), satisfies for any
integer $n\geq1$ the following inequality
$$\mathbb{E}\left[A(\tilde{f}_n)-A^* \right]\leq
\min_{f\in\cC}(A(f)-A^*)+C_0\gamma(\cF,\pi,n,\kappa),$$ where
$C_0>0$ depends only on the constants $\kappa$ and $c_0$ appearing
in MA($\kappa$).

Moreover, there exists a set of prediction rules
$\cF=\{f_1,\ldots,f_M\}$ such that for any procedure $\bar{f}_n$
with values in $\mathbb{R}$, there exists a probability measure
$\pi$ satisfying MA($\kappa$) such that for any integers $M,n$
with $\log M\leq n$ we have
$$\mathbb{E}\left[A(\bar{f}_n)-A^* \right]\geq
\min_{f\in\cC}(A(f)-A^*)+C_0'\gamma(\cF,\pi,n,\kappa),$$where
$C_0'>0$ depends only on the constants $\kappa$ and $c_0$
appearing in MA($\kappa$).
\end{theorem}
The hinge loss is linear on $[-1,1]$, thus, model selection
aggregation or convex aggregation are identical problems if we use
the hinge risk and if we aggregate function with values in
$[-1,1]$. Namely, $\min_{f\in\cF}A(f)=\min_{f\in\cC}A(f).$
Moreover, the result of Theorem \ref{theooracleA} is obtained for
the aggregation of functions with values in $[-1,1]$ and not only
for prediction rules. In fact, only functions with values in
$[-1,1]$ have to be considered when we use the hinge loss since,
for any real valued function $f$, we have
$\max(0,1-y\psi(f(x)))\leq \max(0,1-yf(x))$ for all
$x\in\cX,y\in\{-1,1\}$ where $\psi$ is the projection on $[-1,1]$,
thus, $A(\psi(f))-A^*\leq A(f)-A^*.$ Remark that, under
MA($\kappa$), there exists $c>0$ such that,$
  \mathbb{E}\left[|f(X)-f^*(X)| \right]\leq c\left(A(f)-A^*
  \right)^{1/\kappa}$for all functions $f$ on $\cX$ with values in $[-1,1]$ (cf. \cite{lec:06}) .
The proof of Theorem \ref{theooracleA} is not given here by the
lack of space. It can be found in \cite{lec:06}. Instead, we prove
here the following slightly less general result that we will be
further used to construct adaptive minimax classifiers.
\begin{theorem}\label{oracleA} Let $\kappa\geq1$ and let $\cF=\{f_1,\ldots,f_M\}$ be a finite set
of prediction rules with $M\geq3$. We denote by $\cC$ the convex
hull of $\cF$. We assume that $\pi$ satisfies MA($\kappa$). The
aggregate defined in (\ref{aggregate}) with the exponential
weights (\ref{coef1}) (or (\ref{coefrisk})) satisfies for any
integers $n,M$ and any $a>0$ the following inequality
$$\mathbb{E}\left[A(\tilde{f}_n)-A^* \right]\leq (1+a)
\min_{f\in\cC}(A(f)-A^*)+C\left(\frac{\log
M}{n}\right)^{\frac{\kappa}{2\kappa-1}},$$ where $C>0$ is a
constant depending only on $a$.
\end{theorem}

\begin{corollary}\label{oracleR} Let $\kappa\geq1$, $M\geq3$ and
$\{f_1,\ldots,f_M\}$ be a finite set of prediction rules. We
assume that $\pi$ satisfies MA($\kappa$). The AEW procedure
satisfies for any number $a>0$ and any integers $n,M$ the
following inequality, with $C>0$ a constant depending only on $a$,
$$\mathbb{E}\left[R(\tilde{f}_n)-R^* \right]\leq
2(1+a)\min_{j=1,\ldots,M}(R(f_j)-R^*)+C\left(\frac{\log
M}{n}\right)^{\frac{\kappa}{2\kappa-1}}.$$
\end{corollary}

We denote by $\cP_\kappa$ the set of all probability measures on
$\cX\times\{-1,1\}$ satisfying the margin assumption MA($\kappa$).
Combining Corollary \ref{oracleR} and the following theorem, we
get that the residual
$$\left(\frac{\log M}{n}\right)^{\frac{\kappa}{2\kappa-1}}$$ is
a near optimal rate of model selection aggregation in the sense of
\cite{lec:06} when the underlying probability measure $\pi$
belongs to $\cP_\kappa$.

\begin{theorem}\label{optimalkappa}
For any integers $M$ and $n$ satisfying $M\leq \exp(n)$, there
exists $M$ prediction rules $f_1,\ldots,f_M$ such that for any
classifier $\hat{f}_n$ and any $a>0$, we have
$$\sup_{\pi\in\cP_\kappa}
\left[\mathbb{E}\left[R(\hat{f}_n)-R^*
\right]-2(1+a)\min_{j=1,\ldots,M}(R(f_j)-R^*) \right]\geq C_1
\left(\frac{\log M}{n} \right)^{\frac{\kappa}{2\kappa-1}},$$ where
$C_1=c_0^\kappa/(4e2^{2\kappa(\kappa-1)/(2\kappa-1)}(\log
2)^{\kappa/(2\kappa-1)})$.
\end{theorem}

\par

\section{Adaptivity Both to the Margin and to Regularity.}\label{example}
In this section we give two applications of the oracle inequality
stated in Corollary \ref{oracleR}. First, we construct classifiers
with minimax rates of convergence and second, we obtain adaptive
classifiers by aggregating the minimax ones. Following
\cite{at:05}, we focus on the regularity model where $\eta$
belongs to the H{\"o}lder class.

For any multi-index $s=(s_1,\ldots,s_d)\in\mathbb{N}^d$ and any
$x=(x_1,\ldots,x_d)\in\mathbb{R}^d$, we define $|s|=\sum_{j=1}^d
s_i, s!=s_1!\ldots s_d!, x^s=x_1^{s_1}\ldots x_d^{s_d}$ and
$||x||=(x_1^2+\ldots+x_d^2)^{1/2}.$ We denote by $D^s$ the
differential operator $\frac{\partial^{s_1+\ldots+s_d}}{\partial
x_1^{s_1}\ldots\partial x_d^{s_d}}.$

Let $\beta>0$. We denote by $\lfloor \beta \rfloor$ the maximal
integer that is strictly less than $\beta.$ For any $x\in(0,1)^d$
and any $\lfloor \beta \rfloor$-times continuously differentiable
real valued function $g$ on $(0,1)^d,$ we denote by $g_x$ its
Taylor polynomial of degree $\lfloor \beta \rfloor$ at point $x$,
namely,
$$g_x(y)=\sum_{|s|\leq \lfloor \beta
\rfloor}\frac{(y-x)^s}{s!}D^s g(x).$$

For all $L>0$ and $\beta>0$. The $(\beta,L,[0,1]^d)-${\it
H{\"o}lder class} of functions, denoted by
$\Sigma(\beta,L,[0,1]^d)$, is the set of all real valued functions
$g$ on $[0,1]^d$ that are $\lfloor \beta \rfloor$-times
continuously differentiable on $(0,1)^d$ and satisfy, for any
$x,y\in(0,1)^d,$ the inequality
$$|g(y)-g_x(y)|\leq L||x-y||^\beta.$$

A control of the complexity of H{\"o}lder classes is given by
Kolmogorov and Tikhomorov (1961):
\begin{equation}\label{entropy}
\cN\left(\Sigma(\beta,L,[0,1]^d),\epsilon,L^{\infty}([0,1]^d)
\right)\leq A(\beta,d)\epsilon^{-\frac{d}{\beta}}, \forall
\epsilon>0,\end{equation} where the LHS is the $\epsilon-$entropy
of the $(\beta,L,[0,1]^d)-$H{\"o}lder class w.r.t. to the
$L^\infty([0,1]^d)-$norm and $A(\beta,d)$ is a constant depending
only on $\beta$ and $d$.

If we want to use entropy assumptions on the set which $\eta$
belongs to, we need to make a link between $P^X$ and the Lebesgue
measure, since the distance in (\ref{entropy}) is the
$L^\infty-$norm w.r.t. the Lebesgue measure. Therefore, introduce
the following assumption:

\noindent {\bf (A1)}{\it  The marginal distribution $P^X$ on $\cX$
of $\pi$ is absolutely continuous w.r.t. the Lebesgue measure
$\lambda_d$ on $[0,1]^d$, and there exists a version of its
density which is upper bounded by $\mu_{max}<\infty$. }

We consider the following class of models. For all $\kappa\geq1$
and $\beta>0$, we denote by $\cP_{\kappa,\beta},$ the set of all
probability measures $\pi$ on $\cX\times\{-1,1\}$, such that
\begin{enumerate}\item MA($\kappa$)
is satisfied.\item The marginal $P^X$ satisfies (A1).\item The
conditional probability function $\eta$ belongs to
$\Sigma(\beta,L,\mathbb{R}^d)$.\end{enumerate}

Now, we define the class of classifiers which attain the optimal
rate of convergence, in a minimax sense, over the models
$\cP_{\kappa,\beta}.$ Let $\kappa\geq 1$ and $\beta>0$. For any
$\epsilon>0$, we denote by $\Sigma_\epsilon(\beta)$ an
$\epsilon$-net on $\Sigma(\beta,L,[0,1]^d)$ for the
$L^\infty-$norm, such that, its cardinal satisfies $\log {\rm
Card}\left(\Sigma_\epsilon(\beta) \right)\leq
A(\beta,d)\epsilon^{-d/\beta}$. We consider the AEW procedure
defined in (\ref{aggregate}), over the net
$\Sigma_\epsilon(\beta):$
\begin{equation}\label{classifier}\tf^{\epsilon}=\sum_{\eta\in\Sigma_\epsilon(\beta)}w^{(n)}(f_\eta)f_\eta,
\mbox{ where }
f_\eta(x)=2\1_{\left(\eta(x)\geq1/2\right)}-1.\end{equation}

\begin{theorem}\label{ratecv}
  Let $\kappa>1$ and $\beta>0$. Let $a_1>0$ be an absolute
  constant and consider
  $\epsilon_n=a_1n^{-\frac{\beta(\kappa-1)}{\beta(2\kappa-1)+d(\kappa-1)}}.$
  The aggregate (\ref{classifier}) with $\epsilon=\epsilon_n$,
  satisfies, for any $\pi\in\cP_{\kappa,\beta}$ and any integer $n\geq1$, the following inequality $$
  \mathbb{E}_\pi\left[R(\tf^{\epsilon_n})-R^* \right]\leq C_2(\kappa,\beta,d)
  n^{-\frac{\beta\kappa}{\beta(2\kappa-1)+d(\kappa-1)}},$$ where
  $C_2(\kappa,\beta,d)=2\max\left(4(2c_0\mu_{max})^{\kappa/(\kappa-1)},
  CA(\beta,d)^{\frac{\kappa}{2\kappa-1}}
  \right)(a_1)^{\frac{\kappa}{\kappa-1}}\vee(a_1)^{-\frac{d\kappa}{\beta(\kappa-1)}}$
  and $C$ is the constant appearing in Corollary \ref{oracleR}.
\end{theorem}
Audibert and Tsybakov (cf. \cite{at:05}) have shown the
optimality, in a minimax sense, of the rate obtained in theorem
\ref{ratecv}. Note that this rate is a fast rate because it can
approach $1/n$ when $\kappa$ is close to $1$ and $\beta$ is large.

The construction of the classifier $\tf^{\epsilon_n}$ needs the
knowledge of $\kappa$ and $\beta$ which are not available in
practice. Thus, we need to construct classifiers independent of
these parameters and which learn with the optimal rate
$n^{-\beta\kappa/(\beta(2\kappa-1)+d(\kappa-1))}$ if the
underlying probability measure $\pi$ belongs to
$\cP_{\kappa,\beta}$, for different values of $\kappa$ and
$\beta$. We now show that using the procedure (\ref{aggregate}) to
aggregate the classifiers $\tf^{\epsilon}$, for different values
of $\epsilon$ in a grid, the oracle inequality of Corollary
\ref{oracleR} provides the result.

We use a split of the sample for the adaptation step. Denote by
$D_m^{(1)}$ the subsample containing the first $m$ observations
and $D_l^{(2)}$ the one containing the $l$($=n-m$) last ones.
Subsample $D_m^{(1)}$ is used to construct the classifiers
$\tilde{f}_m^{\epsilon}$ for different values of $\epsilon$ in a
finite grid. Subsample $D_l^{(2)}$ is used to aggregate these
classifiers by the procedure (\ref{aggregate}). We take
$$l=\left\lceil \frac{n}{\log n}\right\rceil \quad \mbox{ and } \quad m=n-l.$$
Set $\Delta=\log n$. We consider a grid of values for $\epsilon$:
$$\cG(n)=\left\{\phi_{n,k}=\frac{k}{\Delta}:k\in\left\{1,\ldots,\lfloor \Delta/2 \rfloor\right\} \right\}.$$
For any $\phi\in\cG(n)$ we consider the step
$\epsilon_m^{(\phi)}=m^{-\phi}.$ The classifier that we propose is
the sign of
$$\tilde{f}_n^{adp}=\sum_{\phi\in\cG(n)}w^{[l]}(\tilde{F}_m^{\epsilon_m^{(\phi)}})\tilde{F}_m^{\epsilon_m^{(\phi)}},$$
where $\tilde{F}_m^{\epsilon}(x)={\rm
sign}(\tilde{f}_m^{\epsilon}(x))$ is the classifier associated to
the aggregate $\tilde{f}_m^{\epsilon}$ for all $\epsilon>0$ and
the weights $w^{[l]}(F)$ are the ones introduced in (\ref{coef1})
constructed with the observations $D_l^{(2)}$ for all
$F\in\cF(n)=\{{\rm
sign}(\tilde{f}_m^{\epsilon}):\epsilon=m^{-\phi},
  \phi\in\cG(n)\}$:
$$w^{[l]}(F)=\frac{\exp\left(
\sum_{i=m+1}^nY_iF(X_i)\right)}{\sum_{G\in\cF(n)}\exp\left(
\sum_{i=m+1}^nY_iG(X_i)\right)}.$$ The following Theorem shows
that $\tilde{f}_n^{adp}$ is adaptive both to the low noise
exponent $\kappa$ and to the complexity (or regularity) parameter
$\beta$, provided that $(\kappa,\beta)$ belongs to a compact
subset of $(1,+\infty)\times(0,+\infty).$

\begin{theorem}\label{adaptation}
  Let $K$ be a compact subset of $(1,+\infty)\times(0,+\infty)$.
  There exists a constant $C_3>0$ that depends
  only on $K$ and $d$ such that for any integer $n\geq1$,
  any $(\kappa,\beta)\in K$ and any
  $\pi\in\cP_{\kappa,\beta}$, we have,
  $$\mathbb{E}_\pi\left[R(\tilde{f}_n^{adp})-R^* \right]\leq C_3n^{-\frac{\kappa\beta}{\beta(2\kappa-1)+d(\kappa-1)}}.$$
\end{theorem}

Classifiers $\tf^{\epsilon_n}$ are not easily implementable since
the cardinality of $\Sigma_{\epsilon_n}(\beta)$ is an exponential
of $n$. An alternative procedure which is easily implementable is
to aggregate plug-in classifiers constructed in Audibert and
Tsybakov (cf. \cite{at:05}).

We introduce the class of models $\cP'_{\kappa,\beta}$ composed of
all the underlying probability measures $\pi$ such that:
\begin{enumerate}
  \item $\pi$ satisfies the margin assumption
        MA($\kappa$).
  \item The conditional probability function
        $\eta\in\Sigma(\beta,L,[0,1]^d).$
  \item  The marginal distribution of $X$ is supported on $[0,1]^d$ and has a Lebesgue density lower bounded  and upper bounded
        by two constants.
\end{enumerate}
\begin{theorem}[{\bf Audibert and Tsybakov (2005)}]\label{theoAT05}
Let $ \kappa>1, \beta>0$. The excess risk of the plug-in
classifier
$\hat{f}^{(\beta)}_n=2\1_{\{\hat{\eta}_n^{(\beta)}\geq1/2\}}-1$
satisfies
$$\sup_{\pi\in\cP'_{\kappa,\beta}}\mathbb{E}\left[R(\hat{f}^{(\beta)}_n)-R^* \right]
\leq C_4n^{-\frac{\beta\kappa}{(\kappa-1)(2\beta+d)}},$$ where
$\hat{\eta}_n^{(\beta)}(\cdot)$ is the locally polynomial
estimator of $\eta(\cdot)$ of order $\lfloor \beta \rfloor$ with
bandwidth $h=n^{-\frac{1}{2\beta+d}}$ and $C_4$ a positive
constant.
\end{theorem}
In \cite{at:05}, it is shown that the rate
$n^{-\frac{\beta\kappa}{(\kappa-1)(2\beta+d)}}$ is minimax over
$\cP'_{\kappa,\beta}$, if $\beta\leq d(\kappa-1)$. Remark that the
fast rate $n^{-1}$ can be achieved.

We aggregate the classifiers $\hat{f}_n^{(\beta)}$ for different
values of $\beta$ lying in a finite grid.  We use a split of the
sample to construct our adaptive classifier: $l=\left\lceil n/\log
n\right\rceil\mbox{ and } m=n-l.$ The training sample
$D_m^{1}=\left((X_1,Y_1),\ldots,(X_m,Y_m)\right)$ is used for the
construction of the class of plug-in classifiers
$$\cF=\left\{\hat{f}_m^{(\beta_k)}:
\beta_k=\frac{kd}{\Delta-2k}, k\in\left\{1,\ldots,\lfloor \Delta/2
\rfloor\right\}\right\},\, \mbox{ where } \Delta=\log n.$$ The
validation sample
$D_l^{2}=\left((X_{m+1},Y_{m+1}),\ldots,(X_n,Y_n)
  \right) \mbox{ }$
is used for the construction of weights
  $$w^{[l]}(f)=\frac{\exp\left(
\sum_{i=m+1}^nY_if(X_i)\right)}{\sum_{\bar{f}\in\cF}\exp\left(
\sum_{i=m+1}^nY_i\bar{f}(X_i)\right)}, \quad \forall f\in\cF.$$
The classifier that we propose is $\tilde{F}_n^{adp}={\rm
sign}(\tilde{f}_n^{adp})$, where:
$\tilde{f}_n^{adp}=\sum_{f\in\cF}w^{[l]}(f)f.$
\begin{theorem}\label{theoadap2}
  Let $K$ be a compact subset of $(1,+\infty)\times(0,+\infty)$.
  There exists a constant $C_5>0$ depending
  only on $K$ and $d$ such that for any integer $n\geq1$,
  any $(\kappa,\beta)\in K$, such that $\beta< d(\kappa-1)$, and any
  $\pi\in\cP'_{\kappa,\beta}$, we have,
  $$\mathbb{E}_\pi\left[R(\tilde{F}_n^{adp})-R^* \right]\leq C_5n^{-\frac{\beta\kappa}{(\kappa-1)(2\beta+d)}}.$$
\end{theorem}
Adaptive classifiers are obtained in Theorem (\ref{adaptation})
and (\ref{theoadap2}) by aggregation of only $\log n$ classifiers.
Other construction of adaptive classifiers can be found in
\cite{lec:05}. In particular, adaptive SVM classifiers.

\par

\section{Proofs}\label{proofs}
\indent {\bf Proof of Proposition \ref{propAEWERM}.} Using the
convexity of the hinge loss, we have
$A_n(\tilde{f_n})\leq\sum_{j=1}^{M}w_j A_n(f_j)$. Denote by
$\hat{i}= {\rm arg}\min_{i=1,\ldots,M} A_n(f_i)$, we have
$A_n(f_i)=A_n(f_{\hat{i}})+\frac{1}{n} \left(
\log(w_{\hat{i}})-\log (w_i)\right)$ for all $i=1,\ldots,M$ and by
averaging over the $w_i$ we get :
 \begin{equation}
 A_n(\tilde{f_n})\leq
 \min_{i=1,\ldots,M}A_n(f_i)+\frac{\log(M)}{n},
 \end{equation}
where we used that $\sum_{j=1}^Mw_j\log \left(\frac{w_j}{1/M}
\right)=K(w|u) \geq 0$ where $K(w|u)$ denotes the Kullback-Leiber
divergence between the weights $w=(w_j)_{j=1,\ldots,M}$ and
uniform weights $u=(1/M)_{j=1,\ldots,M}$.

{\bf Proof of Theorem \ref{oracleA}.} Let $a>0$. Using Proposition
\ref{propAEWERM}, we have for any $f\in\cF$ and for the Bayes rule
$f^*$:
$$ A(\tilde{f}_n)-A^* =
  (1+a)(A_n(\tilde{f}_n)-A_n(f^*))+A(\tilde{f}_n)-A^*-(1+a)(A_n(\tilde{f}_n)-A_n(f^*))$$
$$ \leq  (1+a)(A_n(f)-A_n(f^*))+(1+a)\frac{\log M}{n}\\
+A(\tilde{f}_n)-A^*-(1+a)(A_n(\tilde{f}_n)-A_n(f^*)).$$ Taking the
expectations, we get
\begin{eqnarray*}
\mathbb{E}\left[A(\tilde{f}_n)-A^* \right] & \leq &
(1+a)\min_{f\in\cF}(A(f)-A^*)+(1+a)(\log M)/n\\ & & +
\mathbb{E}\left[A(\tilde{f}_n)-A^*-(1+a)(A_n(\tilde{f}_n)-A_n(f^*))
\right].
\end{eqnarray*}
The following inequality follows from the linearity of the hinge
loss on $[-1,1]$:
$$A(\tilde{f}_n)-A^*-(1+a)(A_n(\tilde{f}_n)-A_n(f^*))\leq
 \max_{f\in\cF}\left[A(f)-A^*-(1+a)(A_n(f)-A_n(f^*)) \right].$$
Thus, using Bernstein's inequality, we have for all
$0<\delta<4+2a:$
\begin{eqnarray*}
  \lefteqn{\mathbb{P}\left[A(\tilde{f}_n)-A^*-(1+a)(A_n(\tilde{f}_n)-A_n(f^*))\geq \delta
  \right]}\\
  &\leq & \sum_{f\in\cF}\mathbb{P}\left[A(f)-A^*-(A_n(f)-A_n(f^*))\geq
  \frac{\delta+a(A(f)-A^*)}{1+a}
  \right]\\
  & \leq & \sum_{f\in\cF}
  \exp\left(-\frac{n(\delta+a(A(f)-A^*))^2}{2(1+a)^2(A(f)-A^*)^{1/\kappa}+2/3(1+a)(\delta+a(A(f)-A^*))}
  \right).
\end{eqnarray*}
There exists a constant $c_1>0$ depending only on $a$ such that
for all $0<\delta<4+2a$ and all $f\in\cF$, we
have$$\frac{(\delta+a(A(f)-A^*))^2}{2(1+a)^2(A(f)-A^*)^{1/\kappa}+2/3(1+a)(\delta+a(A(f)-A^*))}\geq
c_1 \delta^{2-1/\kappa}.$$ Thus,
  $\mathbb{P}\left[A(\tilde{f}_n)-A^*-(1+a)(A_n(\tilde{f}_n)-A_n(f^*))\geq \delta
  \right]\leq M\exp(-nc_1\delta^{2-1/\kappa}).$

Observe that an integration by parts leads to
$\int_a^{+\infty}\exp\left(-bt^\alpha \right)dt\leq
\frac{\exp(-ba^\alpha)}{\alpha b a^{\alpha-1}}$, for any
$\alpha\geq1$ and $a,b>0$, so for all $u>0$, we get
$$\mathbb{E}\left[ A(\tf)-A^*-(1+a)(A_n(\tf)-A_n(f^*))\right]
\leq 2u +M\frac{\exp(-nc_1u^{2-1/\kappa})}{nc_1u^{1-1/\kappa}}.$$
If we denote by $\mu(M)$ the unique solution of $X=M\exp(-X)$, we
have $\log M/2\leq \mu(M)\leq \log M$. For $u$ such that
$nc_1u^{2-1/\kappa}=\mu(M)$, we obtain the result.

{\bf Proof of Corollary \ref{oracleR}.} We deduce Corollary
\ref{oracleR} from Theorem \ref{oracleA}, using that for any
prediction rule $f$ we have $A(f)-A^*=2(R(f)-R^*)$ and applying
Zhang's inequality $A(g)-A^*\geq(R(g)-R^*)$ fulfilled by all $g$
from $\cX$ to $\mathbb{R}$.

{\bf Proof of Theorem \ref{optimalkappa}.} For all prediction
rules $f_1,\ldots,f_M$, we have
$$\sup_{f_1,\ldots,f_M}\inf_{\hat{f}_n}\sup_{\pi\in\cP_\kappa}
\left(\mathbb{E}\left[R(\hat{f}_n)-R^*
\right]-2(1+a)\min_{j=1,\ldots,M}(R(f_j)-R^*) \right)$$
$$ \geq
\inf_{\hat{f}_n}\sup_{\pi\in\cP_\kappa: f^*\in\{f_1,\ldots,f_M\}}
\left(\mathbb{E}\left[R(\hat{f}_n)-R^* \right]\right).$$ Thus, we
look for a set of cardinality not greater than $M$, of the worst
probability measures $\pi\in\cP_\kappa$ from our classification
problem point of view and choose $f_1,\ldots,f_M$ as the
corresponding Bayes rules.

Let $N$ be an integer such that $2^{N-1}\leq M$. Let
$x_1,\ldots,x_N$ be $N$ distinct points of $\cX$. Let $0<w<1/N$.
Denote by $P^X$ the probability measure on $\cX$ such that
$P^X(\{x_j\})=w$ for $j=1,\ldots,N-1$ and $P^X(\{x_N\})=1-(N-1)w$.
We consider the set of binary sequences $\Omega=\{-1,1\}^{N-1}$.
Let $0<h<1$. For all $\sigma\in\Omega$ we consider
$$\eta_{\sigma}(x)=\left\{\begin{array}{ll} (1+\sigma_j h)/2 & \mbox{ if } x=x_1,\ldots,x_{N-1},\\
1 & \mbox{ if } x =x_N. \end{array}
 \right.$$
For all $\sigma\in\Omega$ we denote by $\pi_\sigma$ the
probability measure on $\cX\times\{-1,1\}$ with the marginal $P^X$
on $\cX$ and with  the conditional probability function
$\eta_{\sigma}$ of $Y=1$ knowing $X$.

Assume that $\kappa>1$. We have
$\mathbb{P}\left(|2\eta_\sigma(X)-1|\leq t
\right)=(N-1)w\1_{\{h\leq t\}},\forall 0\leq t<1$. Thus, if we
assume that $(N-1)w\leq h^{1/(\kappa-1)}$ then
$\mathbb{P}\left(|2\eta_\sigma(X)-1|\leq t \right)\leq
t^{1/(\kappa-1)},$ for all $t\geq0$, and according to
\cite{tsy:04}, $\pi_\sigma$ belongs to MA($\kappa$).

We denote by $\rho$ the Hamming distance on $\Omega$ (cf.
\cite{tsybook:04} p.88). Let $\sigma, \sigma'$ be such that
$\rho(\sigma,\sigma')=1$. We have
$$H^2\left(\pi_\sigma^{\otimes n}, \pi_{\sigma'}^{\otimes
n}\right)=2\left(1-(1-w(1-\sqrt{1-h^2}))^n\right).$$ We take $w$
and $h$ such that $w(1-\sqrt{1-h^2})\leq1/n,$ thus,
$H^2\left(\pi_\sigma^{\otimes n}, \pi_{\sigma'}^{\otimes
n}\right)\leq \beta=2(1-e^{-1})<2$ for any integer $n$.

Let $\hat{f}_n$ be a classifier and $\sigma\in\Omega$. Using
MA($\kappa$), we have
$$\mathbb{E}_{\pi_\sigma}\left[ R(\hat{f}_n)-R^*\right]\geq (c_0w)^\kappa\mathbb{E}_{\pi_\sigma}\left[
\left(
\sum_{i=1}^{N-1}|\hat{f}_n(x_i)-\sigma_i|\right)^{\kappa}\right].$$
By Jensen's Lemma and Assouad's Lemma (cf. \cite{tsybook:04}) we
obtain:

$$\inf_{\hat{f}_n}\sup_{\pi\in\cP_\kappa: f^*\in\{f_\sigma:\sigma\in\Omega\}}
\left(\mathbb{E}_{\pi_\sigma}\left[R(\hat{f}_n)-R^*
\right]\right)\geq (c_0w)^\kappa\left(\frac{N-1}{4}(1-\beta/2)^2
\right)^\kappa.$$

We obtain the result by taking $w=(nh^2)^{-1}$, $N=\lceil \log M/
\log 2 \rceil$ and $h=\left(n^{-1}\lceil \log M/ \log 2 \rceil
\right)^{(\kappa-1)/(2\kappa-1)}$.

For $\kappa=1$, we take $h=1/2$, thus $|2\eta_\sigma(X)-1|\geq
1/2$ a.s. so $\pi_\sigma\in$MA(1) (cf.\cite{tsy:04}). Putting
$w=4/n$ and $N=\lceil \log M/ \log 2 \rceil$ we obtain the result.

{\bf Proof of Theorem \ref{ratecv}.} According to Theorem
\ref{oracleR}, where we set $a=1$, we have, for any $\epsilon>0$:
$$\mathbb{E}_\pi\left[R(\tilde{f}_n^{\epsilon})-R^* \right]\leq
4\min_{\bar\eta\in\Sigma_\epsilon(\beta)}\left(R(f_{\bar\eta})-R^*
\right)+C\left(\frac{\log {\rm Card}\Sigma_\epsilon(\beta)}{n}
\right)^{\frac{\kappa}{2\kappa-1}}.$$

Let $\bar\eta$ be a function with values in $[0,1]$ and denote by
$\bar{f}=\1_{\bar{\eta}\geq1/2}$ the plug-in classifier
associated. We have $|2\eta-1|\1_{\bar{f}\neq f^*}\leq
2|\bar{\eta}-\eta|$, thus:
$$R(\bar{f})-R^*  =  \mathbb{E}\left[ |2\eta(X)-1|\1_{\bar{f}\neq
f^*}\right] = \mathbb{E}\left[ |2\eta(X)-1|\1_{\bar{f}\neq
f^*}\1_{\bar{f}\neq f^*}\right]$$ $$\leq  \left\vert\left\vert
|2\eta-1|\1_{\bar{f}\neq f^*}
\right\vert\right\vert_{L^\infty(P^X)}\mathbb{E}\left[
\1_{\bar{f}\neq f^*}\right] \leq   \left\vert\left\vert
|2\eta-1|\1_{\bar{f}\neq f^*}
\right\vert\right\vert_{L^\infty(P^X)} c_0 \left(
R(\bar{f})-R^*\right)^{\frac{1}{\kappa}}, $$
and assumption (A1) lead to
$$R(f_{\bar\eta})-R^*\leq (2c_0\mu_{max})^{\frac{\kappa}{\kappa-1}}||\bar\eta-\eta
||_{L^\infty([0,1]^d)}^{\frac{\kappa}{\kappa-1}}.$$ Hence, for any
$\epsilon>0$, we have
$$\mathbb{E}_\pi\left[R(\tilde{f}_n^{\epsilon})-R^* \right]\leq
 D\left(\epsilon^{\frac{\kappa}{\kappa-1}}+
 \left(\frac{\epsilon^{-d/\beta}}{n}\right)^{\frac{\kappa}{2\kappa-1}}\right),$$
 where $D=\max\left(4(2c_0\mu_{max})^{\kappa/(\kappa-1)}, CA(\beta,d)^{\frac{\kappa}{2\kappa-1}}
 \right)$. For the value $$\epsilon_n=a_1
 n^{-\frac{\beta(\kappa-1)}{\beta(2\kappa-1)+d(\kappa-1)}},$$we have
$$\mathbb{E}_\pi\left[R(\tilde{f}_n^{\epsilon_n})-R^* \right]\leq
C_1 n^{-\frac{\beta\kappa}{\beta(2\kappa-1)+d(\kappa-1)}},$$ where
$C_1=2D(a_1)^{\frac{\kappa}{\kappa-1}}\vee(a_1)^{-\frac{d\kappa}{\beta(\kappa-1)}}$

{\bf Proof of Theorem \ref{adaptation}.} We consider the following
function on $(1,+\infty)\times(0,+\infty)$ with values in
$(0,1/2)$:
$$\phi(\kappa,\beta)=\frac{\beta(\kappa-1)}{\beta(2\kappa-1)+d(\kappa-1)}.$$
For any $n$ greater than $n_1=n_1(K)$, we have $\Delta^{-1}\leq
\phi(\kappa,\beta) \leq \left\lfloor \Delta/2
\right\rfloor\Delta^{-1}$ for all $(\kappa,\beta)\in K.$

Let $(\kappa_0,\beta_0)\in K$. For any $n\geq n_1$, there exists
$k_0\in\{1,\ldots,\lfloor \Delta/2 \rfloor-1\}$ such that
$$\phi_{k_0}=k_0\Delta^{-1}\leq
\phi(\kappa_0,\beta_0)<(k_0+1)\Delta^{-1}.$$ We denote by
$f_{\kappa_0}(\cdot)$ the increasing function
$\phi(\kappa_0,\cdot)$ from $(0,+\infty)$ to $(0,1/2)$. We set
$$\beta_{0,n}=\left(f_{\kappa_0}\right)^{-1}(\phi_{k_0}).$$
There exists $m=m(K)$ such that
$m|\beta_0-\beta_{0,n}|\leq|f_{\kappa_0}(\beta_0)-f_{\kappa_0}(\beta_{0,n})|\leq
\Delta^{-1}.$

Let $\pi\in\cP_{\kappa_0,\beta_0}.$ According to the oracle
inequality of Corollary \ref{oracleR}, we have, conditionally to
the first subsample $D_m^{1}$:
$$\mathbb{E}_\pi\left[R(\tilde{f}_n^{adp})-R^*| D_m^{1} \right]\leq
 4
 \min_{\phi\in\cG(n)}\left(R(\tilde{f}_m^{\epsilon_m^{(\phi)}})-R^*\right)
 +C\left(\frac{\log {\rm Card}(\cG(n))}{l}
 \right)^{\frac{\kappa_0}{2\kappa_0-1}}.$$
Using the definition of $l$ and the fact that ${\rm
Card}(\cG(n))\leq\log n$ we get that there exists $\tilde{C}>0$
independent of $n$ such that
 $$\mathbb{E}_\pi\left[R(\tilde{f}_n^{adp})-R^*\right]\leq\tilde{C}\left(
 \mathbb{E}_\pi\left[R(\tilde{f}_m^{\epsilon_m^{(\phi_{k_0})}})-R^* \right]+
 \left(\frac{\log^2 n}{n} \right)^{\frac{\kappa_0}{2\kappa_0-1}}\right)$$

Moreover $\beta_{0,n}\leq \beta_0$, hence,
$\cP_{\kappa_0,\beta_0}\subseteq \cP_{\kappa_0,\beta_{0,n}}$.
Thus, according to Theorem \ref{ratecv}, we have
$$\mathbb{E}_\pi\left[R(\tilde{f}_m^{\epsilon_m^{(\phi_{k_0})}})-R^*
\right]\leq C_1(K,d)m^{-\psi(\kappa_0,\beta_{0,n})},$$ where
$C_1(K,d)=\max\left(C_1(\kappa,\beta,d):(\kappa,\beta)\in K
\right)$ and
$\psi(\kappa,\beta)=\frac{\beta\kappa}{\beta(2\kappa-1)+d(\kappa-1)}.$
By construction, there exists $A_2=A_2(K,d)>0$ such that
$|\psi(\kappa_0,\beta_{0,n})-\psi(\kappa_0,\beta_0)|\leq
A_2\Delta^{-1}.$ Moreover for any integer $n$ we have $n^{A_2/\log
n}=\exp(A_2)$, which is a constant. We conclude that
$$\mathbb{E}_\pi\left[R(\tilde{f}_n^{adp})-R^*
\right]\leq C_2(K,d)
\left(n^{-\psi(\kappa_0,\beta_0)}+\left(\frac{\log^2 n}{n}
\right)^{\frac{\kappa_0}{2\kappa_0-1}} \right),$$ where
$C_2(K,d)>0$ is independent of $n$. We achieve the proof by
observing that
$\psi(\kappa_0,\beta_0)<\frac{\kappa_0}{2\kappa_0-1}.$

{\bf Proof of Theorem \ref{theoadap2}.} We consider the following
function on $(1,+\infty)\times(0,+\infty)$ with values in
$(0,1/2)$:
$$\Theta(\kappa,\beta)=\frac{\beta\kappa}{(\kappa-1)(2\beta+d)}.$$
For any $n$ greater than $n_1=n_1(K)$, we have $\Delta^{-1}\leq
\Theta(\kappa,\beta) \leq \left\lfloor \Delta/2
\right\rfloor\Delta^{-1},$ for all $(\kappa,\beta)\in K.$

Let $(\kappa_0,\beta_0)\in K$ be such that
$\beta_0<(\kappa_0-1)d.$ For any $n\geq n_1$, there exists
$k_0\in\{1,\ldots,\lfloor \Delta/2 \rfloor-1\}$ such that
$k_0\Delta^{-1}\leq \Theta(\kappa_0,\beta_0)<(k_0+1)\Delta^{-1}.$

Let $\pi\in\cP_{\kappa_0,\beta_0}.$ According to the oracle
inequality of Corollary \ref{oracleR}, we have, conditionally to
the first subsample $D_m^{1}$:
$$\mathbb{E}_\pi\left[R(\tilde{F}_n^{adp})-R^*| D_m^{1} \right]\leq
 4 \min_{f\in\cF}(R(f)-R^*)
 +C\left(\frac{\log {\rm Card}(\cF)}{l}
 \right)^{\frac{\kappa_0}{2\kappa_0-1}}.$$
Using the proof of Theorem \ref{adaptation} we get that there
exists $\tilde{C}>0$ independent of $n$ such that
 $$\mathbb{E}_\pi\left[R(\tilde{f}_n^{adp})-R^*\right]\leq\tilde{C}\left(
 \mathbb{E}_\pi\left[R(\hat{f}_m^{(\beta_{k_0})})-R^* \right]+
 \left(\frac{\log^2 n}{n} \right)^{\frac{\kappa_0}{2\kappa_0-1}}\right)$$

Moreover $\beta_{k_0}\leq \beta_0$, hence,
$\cP_{\kappa_0,\beta_0}\subseteq \cP_{\kappa_0,\beta_{k_0}}$.
Thus, according to Theorem \ref{theoAT05}, we have
$$\mathbb{E}_\pi\left[R(\hat{f}_m^{(\beta_{k_0})})-R^*
\right]\leq C_4(K,d)m^{-\Theta(\kappa_0,\beta_{k_0})},$$ where
$C_4(K,d)=\max\left(C_4(\kappa,\beta,d):(\kappa,\beta)\in K
\right)$. We have
$|\Theta(\kappa_0,\beta_{k_0})-\Theta(\kappa_0,\beta_0)|\leq
\Delta^{-1}$ by construction. Moreover $n^{1/\log n}=e$ for any
integer $n$. We conclude that
$$\mathbb{E}_\pi\left[R(\tilde{F}_n^{adp})-R^*
\right]\leq \tilde{C}_4(K,d)
\left(n^{-\Theta(\kappa_0,\beta_0)}+\left(\frac{\log^2 n}{n}
\right)^{\frac{\kappa_0}{2\kappa_0-1}} \right),$$ where
$\tilde{C}_4(K,d)>0$ is independent of $n$. We achieve the proof
by observing that
$\Theta(\kappa_0,\beta_0)<\frac{\kappa_0}{2\kappa_0-1}$, if
$\beta_0<(\kappa_0-1)d.$


\end{document}